\documentclass[12pt]{article}
\usepackage[english]{babel}
\usepackage{amssymb,amsthm}
\usepackage{amsmath}
\usepackage{inputenc}
\usepackage{graphicx}

\newtheorem{thm}{Theorem}[section]
\newtheorem{lem}{Lemma}[section]
\newtheorem{prop}{Proposition}[section]
\newtheorem{defn}{Definition}[section]

\newtheorem{cor}{Corollary}[section]

\begin{document}

\begin{center}
{\large\bf Classification of complex naturally graded
 \\ quasi-filiform
Zinbiel algebras}\\[2mm]
\end{center}

\begin{center}
{\bf J.Q. Adashev \\
{\small\it Institute of Mathematics and Information Technologies,
\\ Uzbek Academy of Sciences, Tashkent, Uzbekistan,  e-mail:
adashevjq@mail.ru}}
\end{center}

\begin{center}
{\bf  A.Kh. Khudoyberdiyev \\
{\small\it Institute of Mathematics and Information Technologies,
\\ Uzbek Academy of Sciences, Tashkent, Uzbekistan, e-mail:
abror\_alg\_85@mail.ru}}
\end{center}

\begin{center}{ \bf B.A. Omirov \footnote{This work is partially supported by DFG project 436 USB
 113/10/0-1 project (Germany).}\\
{\small\it Institute of Mathematics and Information Technologies,
\\ Uzbek Academy of Sciences, Tashkent, Uzbekistan, e-mail:
omirovb@mail.ru}}
\end{center}

\begin{abstract} In this work the description up to isomorphism of complex naturally graded quasi-filiform Zinbiel
algebras is obtained.

{\bf Key words:} Zinbiel algebra, Leibniz algebra, nilpotency,
nul-filiform algebra, filiform and quasi-filiform algebras.

{\bf AMS classification:} {\rm 17A32.}
\end{abstract}
\begin{center}
This work is dedicated to the 60-th anniversary of Professor
Shestakov I.P.
\end{center}

\section{Introduction}

In the present paper we investigate algebras which are Koszul dual
to Leibniz algebras. The Leibniz algebras were introduced in the
work \cite{Lod2} and they present a "non commutative" (to be more
precise, a "non antisymmetric") analogue of Lie algebras. Many
works, including \cite{Lod2}-\cite{O2}, were devoted to the
investigation of cohomological and structural properties of
Leibniz algebras. Ginzburg and Kapranov introduced and studied the
concept of Koszul dual operads \cite{Gunz1}. Following this
concept, it was shown in \cite{Lod1} that the category of dual
algebras to the category of Leibniz algebras is defined by the
identity:
$$
(x\circ y)\circ z=x\circ(y\circ z)+x\circ (z\circ y).
$$

In this paper, dual Leibniz algebras will called Zinbiel algebras
(Zinbiel is obtained from Leibniz written in inverse order). Some
interesting properties of Zinbiel algebras were obtained in
\cite{Ada}, \cite{Dzhu1}, and \cite{Dzhu2}. In particular, the
nilpotency of an arbitrary finite-dimensional complex Zinbiel
algebra was proved in \cite{Dzhu2}, and zero-filiform and filiform
Zinbiel algebras were classified in \cite{Ada}. The classification
of complex Zinbiel algebras up to dimension $4$ is obtained in
works \cite{Dzhu2} and  \cite{O3}. The present paper is devoted to
the investigation of the next stage
--- description of quasi-filiform complex Zinbiel algebras.

Examples of Zinbiel algebras can be found in \cite{Ada},
\cite{Dzhu2} and \cite{Lod1}.

We consider below only complex algebras and, for convenience, we
will omit zero products the algebra's multiplication table.

\section{Preliminaries}

\begin{defn} An algebra $A$ over a field $F$ is called a Zinbiel
algebra if for any $x, y, z \in A$ the identity
\begin{equation}\label{e1}
(x\circ y)\circ z=x\circ(y\circ z)+x\circ (z\circ y)
\end{equation}
holds.
\end{defn}

For an arbitrary Zinbiel algebras define the lower central series
$$A^1=A,\ \ A^{k+1}=A\circ A^k, \ \ k\geq 1.$$

\begin{defn} A Zinbiel algebra $A$ is called nilpotent if there
exists an $s\in N$ such that $A^s=0$. The minimal number $s$
satisfying this property (i.e. $A^{s-1}\neq 0$ and $A^s=0$) is
called nilindex of the algebra $A$.
\end{defn}

It is not difficult to see that nilindex of an arbitrary
$n$-dimensional nilpotent algebra does not exceed the number
$n+1$.

\begin{defn} An $n$-dimensional Zinbiel algebra $A$ is called
zero-filiform if $dimA^i=(n+1)-i,$ $1 \le i \le n+1.$
\end{defn}

Clearly, the definition of a zero-filiform algebra $A$ amounts to
requiring that $A$ has a maximal nilindex.

\begin{thm} \emph{\cite{Ada}} \label{th1} An arbitrary
$n$-dimensional zero-filiform Zinbiel algebra is isomorphic to the
algebra
\begin{equation} \label{e2}
e_i \circ e_j=C_{i+j-1}^j e_{i+j}, \  \mbox{ for }  \ 2\le i+j \le
n
\end{equation}
where symbol $C_s^t$ is a binomial coefficient defined as
$C_s^t=\frac{s!}{t!(s-t)!},$ and $\{e_1, e_2, \ldots, e_n\}$ is a
basis of the algebra.
\end{thm} We denote the algebra from Theorem \ref{th1} as $NF_n.$

It is easy to see that an $n$-dimensional Zinbiel algebra is
one-generated if and only if it is isomorphic to $NF_n.$

\begin{defn} An $n$-dimensional Zinbiel algebra $A$  is called
filiform if $dimA^i=n-i,$ $ 2 \le i \le n.$
\end{defn}

The following theorem gives classification of filiform Zinbiel
algebras.

\begin{thm} \label{th2} Any $n$-dimensional $(n\geq 5)$ filiform
Zinbiel algebra is isomorphic to one of the following three
pairwise non isomorphic algebras:
$$
\begin{array}{lll} F_n^1: e_i\circ e_j=C_{i+j-1}^j e_{i+j}, & 2 \leq i+j
\leq  n-1;& \\
F_n^2: e_i\circ e_j=C_{i+j-1}^j  e_{i+j}, & 2 \leq i+j \leq
 n-1, & e_n\circ e_1=e_{n-1}; \\
F_n^3: e_i\circ e_j=C_{i+j-1}^j  e_{i+j}, & 2 \leq i+j \leq
 n-1, & e_n\circ e_n=e_{n-1}. \end{array}
$$
\end{thm}

Since the direct sum of nilpotent Zinbiel algebras is nilpotent
one, we shall consider only non split algebras.

Summarizing the results of \cite{Ada},  \cite{Dzhu2}, and
\cite{O3}, we give the classification of complex Zinbiel algebras
up to dimension $\leq 4.$

\begin{thm} \label{th2.3} An arbitrary non split Zinbiel algebra is isomorphic to
the following pairwise non isomorphic algebras:
\\ $DimA=1:$ Abelian
\\ $DimA=2:$ $e_1\circ e_1=e_2$
\\ $DimA=3:$
\begin{description}
\item[$Z_3^1:$] $e_1\circ e_1=e_2, \ \ e_1\circ e_2=\frac{1}{2}e_3, \ \ e_2\circ e_1=e_3;$
\item[$Z_3^2:$] $e_1\circ e_2= e_3, \ \ e_2\circ e_1=-e_3;$
\item[$Z_3^3:$] $e_1\circ e_1= e_3, \ \ e_1\circ e_2=e_3, \ \
     e_2\circ e_2=\alpha e_3, \ \ \alpha \in C;$
\item[$Z_3^4:$] $e_1\circ e_1=e_3, \ \ e_1\circ e_2=e_3, \ \ e_2\circ e_1= e_3.$
\end{description}
$DimA=4:$
\begin{description}
    \item[$Z_4^1:$] $e_1\circ e_1=e_2, \ \ e_1\circ e_2=e_3,\ \ e_2\circ
e_1=2e_3,\ \ e_1\circ e_3=e_4, \ \ e_2\circ e_2=3e_4, \ \ e_3\circ
e_1=3e_4;$
    \item[$Z_4^2:$] $e_1\circ e_1=e_3,\ \  e_1\circ e_2=e_4,\ \ e_1\circ
e_3=e_4, \ \ e_3\circ e_1=2e_4;$
    \item[$Z_4^3:$] $e_1\circ e_1=e_3,\ \ e_1\circ e_3=e_4, \ \ e_2\circ
e_2=e_4, \ \ e_3\circ e_1=2e_4;$
    \item[$Z_4^4:$] $e_1\circ e_2=e_3,\ \ e_1\circ e_3=e_4,\ \ e_2\circ e_1=
-e_3;$
    \item[$Z_4^5:$] $e_1\circ e_2=e_3, \ \ e_1\circ e_3=e_4,\ \ e_2\circ e_1=
-e_3,\ \ e_2\circ e_2= e_4;$
    \item[$Z_4^6:$] $e_1\circ e_1=e_4,\ \ e_1\circ e_2=e_3,\ \ e_2\circ e_1=
-e_3,\ \ e_2\circ e_2= -2e_3+e_4;$
    \item[$Z_4^7:$] $e_1\circ e_2=e_3,\ \ e_2\circ e_1=e_4,\ \ e_2\circ e_2=
-e_3;$
    \item[$Z_4^8(\alpha):$] $ e_1\circ e_1=e_3,\ \ e_1\circ e_2=e_4,\ \
e_2\circ e_1= -\alpha e_3, \ \ e_2\circ e_2= -e_4, \quad \alpha
\in C$
    \item[$Z_4^9(\alpha):$] $e_1\circ e_1=e_4,\ \ e_1\circ e_2=\alpha e_4,\ \
e_2\circ e_1= -\alpha e_4,\ \ e_2\circ e_2=e_4, \ \ e_3\circ
e_3=e_4,\quad  \alpha \in C$
    \item[$Z_4^{10}:$] $e_1\circ e_2=e_4,\ \ e_1\circ e_3=e_4,\ \ e_2\circ
e_1= -e_4,\ \ e_2\circ e_2=e_4,\ \ e_3\circ e_1=e_4;$
    \item[$Z_4^{11}:$] $e_1\circ e_1=e_4,\ \ e_1\circ e_2=e_4,\ \ e_2\circ
e_1= -e_4,\ \ e_3\circ e_3=e_4;$
    \item[$Z_4^{12}:$] $e_1\circ e_2=e_3,\ \ e_2\circ e_1=e_4;$
    \item[$Z_4^{13}:$] $e_1\circ e_2=e_3,\ \ e_2\circ e_1= -e_3, \ \ e_2\circ
e_2=e_4;$
    \item[$Z_4^{14}:$] $e_2\circ e_1=e_4,\ \ e_2\circ e_2=e_3;$
    \item[$Z_4^{15}(\alpha):$] $e_1\circ e_2=e_4,\ \ e_2\circ e_2=e_3,\ \
e_2\circ e_1=\frac{1+\alpha} {1-\alpha}e_4, \quad \alpha \in
C\setminus\{1\};$
    \item[$Z_4^{16}:$] $e_1\circ e_2=e_4,\ \ e_2\circ e_1= -e_4,\ \ e_3\circ
e_3=e_4;$
\end{description}
\end{thm}

Let us introduce some definitions and notations.

The set $R(A)=\{ a\in A | \ \  b\circ a=0$ for any $b\in A\}$ is
called the right annihilator of the Zinbiel algebra $A.$

The set  $L(A)=\{ a\in A | \ \  a\circ b=0$ for any $b\in A\}$ is
called the left annihilator of the Zinbiel algebra $A.$

$Z(a,b,c)$ denotes the following polynomial:
$$
Z(a,b,c)=(a\circ b)\circ c- a\circ (b\circ c)-a\circ (c\circ b).
$$

It is obvious that Zinbiel algebras are determined by the identity
$Z(a,b,c)=0.$

\section{Classification of naturally graded quasi-filiform \\Zinbiel
algebras}

\begin{defn} A Zinbiel algebra $A$ is called quasi-filiform if
$A^{n-2}\ne 0$ and $A^{n-1}= 0,$ where $dimA=n.$
\end{defn}

Let $A$ be a quasi-filiform Zinbiel algebra. Putting
$A_i=A^i/A^{i+1},$ $1 \leq i \leq n-2,$ we obtain the graded
Zinbiel algebra
$$
GrA= A_1\oplus A_2 \oplus \ldots \oplus A_{n-2},   \ \mbox{ where
} \ A_i \circ A_j \subseteq A_{i+j}.
$$

An algebra $A$ is called naturally graded if $A\cong GrA.$ It is
not difficult to see that $A_1\circ A_j=A_{j+1}$ in the naturally
graded algebra $A.$ Let $A$ be an $n$-dimensional graded
quasi-filiform algebra. Then there exists a basis $\{e_1, e_2,
\ldots, e_n\}$ of the algebra $A$ such that $e_i\in A_i,$ $1\le
i\le n-2.$ It is evident that $dimA_1>1.$ In fact, if $dimA_1=1,$
then the algebra $A$ is one-degenerated and therefore it is the
zero-filiform algebra, but it is not quasi-filiform. Without loss
of generality, one can assume $e_{n-1}\in A_1.$ If for a Zinbiel
algebra $A,$ the condition $e_n\in A_r$ holds, the algebra is said
to be of type $A_{(r)}.$

\subsection{The case $r=1$}

\begin{thm} \label{th3} Any $n$-dimensional ($n\geq 6$) naturally
graded quasi-filiform Zinbiel algebra of the type $A_{(1)}$ is
isomorphic to the algebra:
\begin{equation}\label{e3}
KF_n^1: e_i \circ e_j=C_{i+j-1}^j e_{i+j}  \ \mbox{ for }  \ 2\leq
i+j \leq n-2.
\end{equation}
\end{thm}
\begin{proof} Let an algebra $A$ satisfy the conditions of the
theorem. Then there exists a basis $\{e_1, e_2, \ldots, e_n\}$
such that
$$
A_1=\langle e_1, e_{n-1}, e_n\rangle, \ \ A_2=\langle e_2\rangle,
\ \ A_3=\langle e_3\rangle, \ \
 \ldots, \ \ A_{n-2}=\langle e_{n-2}\rangle.
$$
Using arguments similar to the ones from \cite{Ver}, we obtain
$$e_1\circ e_i=e_{i+1} \ \mbox{ for }  \ 2\leq i\leq n-3.$$
Now introduce the notations
$$
\begin{array}{lll}
  e_1\circ e_1=\alpha_{1,1}e_2, & e_1\circ e_{n-1}=\alpha_{1,2} e_2, & e_1\circ e_n=\alpha_{1,3} e_2, \\
  e_{n-1}\circ e_1=\alpha_{2,1}e_2, & e_{n-1}\circ
e_{n-1}=\alpha_{2,2} e_2, & e_{n-1}\circ e_n=\alpha_{2,3} e_2, \\
  e_n\circ e_1=\alpha_{3,1}e_2, & e_n\circ e_{n-1}=\alpha_{3,2}
e_2, & e_n\circ e_n=\alpha_{3,3} e_2. \\
\end{array}$$
We consider the following cases:

\textbf{Case 1.} Let $(\alpha_{1,1},\alpha_{2,2}, \alpha_{3,3})
\neq (0,0,0).$ Then without loss of generality, one can take
$\alpha_{1,1} \neq 0$. Moreover, making the change
$e'_2=\alpha_{1,1}e_2, \ \ e'_3=\alpha_{1,1}e_3, \ \ \ldots, \ \
e'_{n-2}=\alpha_{1,1}e_{n-2},$ we can assume $\alpha_{1,1}=1.$

Obviously, the linear span $lin\langle e_1, e_2, \ldots,
e_{n-2}\rangle$ forms zero-filiform Zinbiel algebra. Hence
$$e_i\circ e_j=  C_{i+j-1}^j e_{i+j} \  \mbox{ at }  \ 2\leq i+j \leq n-2$$
and omitted products of the basic elements $\{e_1, e_2, \ldots,
e_{n-2}\}$ are equal to zero. \\ Taking into account equalities
$$Z(e_1,e_{n-1},e_1)=Z(e_1,e_1,e_{n-1})=Z(e_1,e_{n-1},e_{n-1})=Z(e_{n-1},e_1,e_1)=0$$
we obtain
$$
\alpha_{1,2}=\alpha_{2,1}, \ \ e_2\circ e_{n-1}=2\alpha_{1,2}e_3,
\ \ \alpha_{1,2}^2=\alpha_{2,2}, \ \ e_{n-1}\circ
e_2=\alpha_{1,2}e_3.
$$
Substituting $e'_{n-1}=-\alpha_{1,2}e_1+e_{n-1},$ $e'_i=e_i$ for
$i\neq n-1,$ one can suppose $$\alpha_{1,2}=\alpha_{2,1}
=\alpha_{2,2}=0.$$ Analogously, we have
$$e_1\circ e_n = e_n\circ e_1 = e_n\circ e_n = e_2\circ e_n = e_n\circ e_2 = 0
 \quad \hbox{\rm and} \quad e_n\circ e_{n-1} = e_{n-1}\circ e_n =
0.$$ Hence, using induction and the following chain of equalities
$$e_{n-1}\circ e_i = e_{n-1}\circ (e_1\circ
e_{i-1}) = (e_{n-1}\circ e_1)\circ e_{i-1}-e_{n-1}\circ
(e_{i-1}\circ e_1) =-(i-1)e_{n-1}\circ e_i$$ we prove validity of
the equality $e_{n-1}\circ e_i=0$ for $1\leq i\leq n.$ \\ In the
same way, equalities $$e_{i+1}\circ e_{n-1} = (e_1\circ e_i)\circ
e_{n-1} = e_1\circ (e_i\circ e_{n-1}) + e_1\circ (e_{n-1}\circ
e_i) = 0$$ lead to $e_{i+1}\circ e_{n-1} = 0,$ $1\leq i\leq n.$
But it means that the algebra $A$ is split, i.e. $A =
F_{n-1}^1\oplus C,$ where $F_{n-1}^1 = NF_{n-2}\oplus C.$ Hence,
$A =KF_n^1.$ \\ \textbf{Case 2.} Let
$(\alpha_{1,1},\alpha_{2,2},\alpha_{3,3}) = (0,0,0).$ Then
$$(\alpha_{1,2},\alpha_{2,1},\alpha_{1,3},\alpha_{3,
1},\alpha_{2,3},\alpha_{3,2})\neq (0, 0, 0, 0, 0, 0).$$ Put
$e'_1=ae_1+be_{n-1}+ce_n.$ Then
$$
e'_1\circ e'_1 =
[ab(\alpha_{1,2}+\alpha_{2,1})+ac(\alpha_{1,3}+\alpha_{3,1})+bc
(\alpha_{2,3}+\alpha_{3,2})]e_2.
$$
From this it follows $\alpha_{1,2}+\alpha_{2,1}=0,$ $\alpha_{1,3}+
\alpha_{3,1}=0,$ and $\alpha_{2,3}+\alpha_{3,2}=0.$ In fact,
otherwise revert to the conditions of the case 1. \\Without loss
of generality, we can assume $\alpha_{1,2}=1.$ The equality
$Z(e_1,e_{n-1},e_1)=0$ implies $e_2\circ e_1 =0.$ In addition, $$
0 = (e_1\circ e_1)\circ e_2 = e_1\circ (e_1\circ e_2) + e_1\circ
(e_2\circ e_1) = e_1\circ e_3 = e_4.$$ Thus we obtain a
contradiction with existence of an algebra in case 2. \end{proof}

\begin{prop} Let $A$ be a five-dimensional naturally graded
quasi-filiform Zinbiel algebra of type $A_{(1)}.$ Then it is
isomorphic to one of the following three pairwise non isomorphic
algebras:
$$\begin{array}{l}
KF_5^1: e_1\circ e_1=e_2, \ \ e_1\circ e_2=e_3, \ \ e_2\circ
e_1=e_3; \\[1mm]
KF_5^2: e_1\circ e_4=e_2, \ \ e_1\circ e_2=e_3, \ \ e_4\circ
e_1=-e_3; \\[1mm]
KF_5^3: e_1\circ e_4=e_2, \ \ e_1\circ e_2=e_3, \ \ e_4\circ
e_1=-e_3, \ \ e_4\circ e_1=e_3.\end{array}$$
\end{prop}
\begin{proof} Let an algebra $A$ satisfy the conditions of the proposition.
If the conditions of the case 1 of Theorem \ref{th3} hold, then
$A$ is isomorphic to the algebra $KF_5^1.$ if the conditions of
the case 2 of Theorem \ref{th3} are valid for $A,$ we obtain the
following multiplication in $A:$
$$
\begin{array}{lll} e_1\circ e_4=e_2, & e_1\circ e_5=\alpha_{1,3}e_2, &
e_1\circ e_2=e_3, \\
e_4\circ e_1=-e_2, & e_4\circ e_5=\alpha_{2,3}e_2, & e_4\circ
e_2=\beta_1 e_3, \\
e_5\circ e_1=-\alpha_{1,3}e_2, & e_5\circ e_4=-\alpha_{2,3}e_2, &
e_5\circ e_2=\beta_2e_3.\end{array}
$$

Changing basic elements by the rules
$$e'_1=e_1, \quad e'_2=e_2, \quad e'_3=e_3, \quad
e'_4=-\beta_1e_1+e_4, \quad
e'_5=\alpha_{2,3}e_1-\alpha_{1,3}e_4+e_5,$$ we can assume
$\alpha_{1,3}=\alpha_{2,3}=\beta_1=0.$ If $\beta_2=0,$ then we
obtain the algebra $KF_5^2.$ But if $\beta_2\neq 0,$ then putting
$e'_5= \frac{1}{\beta_2}e_5,$ we have $\beta_2=1$ and obtain
isomorphism to algebra $KF_5^3.$ By virtue of $dimL(KF_5^2)=3,$
$dimL(KF_5^3)=2,$ and taking into account that generating elements
of algebras $KF_5^2$ and $KF_5^3$ satisfy the identity $x\circ
x=0,$ but this identity does not hold for generating elements of
the algebra $KF_5^1$ (in particular, $e_1\circ e_1=e_2$), we
obtain pairwise non isomorphic of obtained algebras.\end{proof}

\subsection{The case $r=2$}

Classification of naturally graded quasi-filiform Zinbiel algebras
of type $A_{(2)}$ is given in the following theorem.

\begin{thm} \label{th4} Any $n$-dimensional $(n\geq 8)$ naturally
graded quasi-filiform Zinbiel algebra of type $A_{(2)}$ is
isomorphic to one of the following pairwise non isomorphic
algebras:
\begin{description}
    \item[$KF^1_n:$] $\left\{\begin{array}{l} e_i\circ e_j=C_{i+j-1}^j e_{i+j}, \ \ 2 \leq i +j \leq n-2,
\\[1mm]
    e_1\circ e_{n-1}=e_n, \ \ e_{n-1}\circ e_1= \alpha e_n;\end{array}\right.$
    \item[$KF^2_n:$] $\left\{\begin{array}{l} e_i\circ e_j=C_{i+j-1}^j
e_{i+j}, \ \ 2 \leq i+j \leq n-2,\\[1mm]
e_1\circ e_{n-1}=e_n, \ \ e_{n-1}\circ e_1=e_n, \ \ e_{n-1}\circ
e_{n-1}=e_n;\end{array}\right.$
    \item[$KF^3_n:$] $\left\{\begin{array}{l} e_i\circ e_j= C_{i+j-1}^j e_{i+j}, \ \
    2 \leq i+j \leq n-2,\\[1mm]
    e_1\circ e_{n-1}=e_n, \ \ e_{n-1}\circ e_{n-1}=e_n;\end{array}\right.$
    \item[$KF^4_n:$] $\left\{\begin{array}{l} e_i\circ e_j = C_{i+j-1}^j
e_{i+j}, \ \ 2 \leq i+j \leq n-2,\\[1mm]
e_{n-1}\circ e_1=e_n.\end{array}\right.$
\end{description}
\end{thm}
\begin{proof} Let an algebra $A$ satisfy the conditions of theorem
and let $\{e_1, e_2, \ldots, e_n\}$ be a basis of $A$ such that
$A_1=\langle e_1, e_{n-1}\rangle,$ $A_2=\langle e_2, e_n\rangle,$
$A_i=\langle e_i\rangle$ for $3\leq i \leq n-2.$ Analogously as in
the proof of Theorem \ref{th3}, we obtain
$$\begin{array}{ll} e_1\circ e_1 = \alpha_1 e_2 + \alpha_2e_n, & e_1\circ e_{n-1}
= \alpha_3e_2 + \alpha_4 e_n, \\[1mm]
e_{n-1}\circ e_1 = \alpha_5e_2 +\alpha_6e_n, & e_{n-1}\circ
e_{n-1} = \alpha_7e_2 + \alpha_8e_n.\end{array}$$
$$\begin{array}{l} e_1\circ e_i =
e_{i+1}  \ \mbox{ for }  \ 3\leq i\leq n-3.\end{array}$$

Without loss of generality one can suppose $e_1\circ e_1 = e_2.$
In fact, if there exists $x \in A_1$ such that $x\circ x \neq 0,$
we can set $e_1=x$ and $e_2=x\circ x.$ But if for any $x\in A_1,$
we have $x\circ x=0,$ then $e_1\circ e_1=e_{n-1}\circ e_{n-1}=0, \
\ e_1\circ e_{n-1}=-e_{n-1}\circ e_1$ and in this case $dim
A_2=1,$ what contradicts to conditions of the theorem.

Thus, $$\begin{array}{ll} e_1\circ e_1 = e_2, & e_1\circ
e_{n-1} = \alpha_3e_2 + \alpha_4e_n, \\[1mm]
 e_{n-1}\circ e_1 =\alpha_5e_2 + \alpha_6e_n, & e_{n-1}\circ e_{n-1} = \alpha_7e_2
+ \alpha_8e_n.\end{array}$$ \\ \textbf{Case 1.} Let $\alpha_4\neq
0$. Then we can assume $e_1\circ e_{n-1} = e_n.$ Let us introduce
denotations  $$e_1\circ e_2=\beta_1e_3, \ \ e_{n-1}\circ
e_2=\beta_2e_3, \ \ e_1\circ e_n=\beta_3e_3, \ \ e_{n-1}\circ
e_n=\beta_4e_5.$$ \\ \textbf{Case 1.1.} \label{Case 1.1} Let
$\beta_1\neq 0.$ Then we can assume $e_1\circ e_2=e_3.$ Putting
$e'_{n-1}=-\beta_2e_1+ e_{n-1},$ we have $\beta_2=0,$ i.e.
$e_{n-1}\circ e_2=0.$ Using the equality $e_1\circ e_i=e_{i+1}$
for $1 \leq i \leq n-3$ and induction, we can get the following
equality
$$e_i\circ e_j = C_{i+j-1}^je_{i+j}  \ \mbox{ for }  \ 2\leq i+j\leq n-2.$$
The equalities
$$\begin{array}{l}
Z(e_1,e_{n-1},e_1)=Z(e_1,e_1,e_{n-1})=Z(e_n,e_1,e_1)=Z(e_1,e_n,e_1)=
\\[1mm]
Z(e_1,e_1,e_n)=Z(e_1,e_{n-1},e_2)=Z(e_{n-1},e_2,e_1)=Z(e_{n-1},e_1,e_2)=
\\[1mm]
Z(e_1,e_{n-1},e_{n-1})=Z(e_{n-1},e_{n-1},e_2)=Z(e_1,e_n,e_{n-1})=Z(e_1,e_{n-1},e_n)=0,\end{array}$$
deduce
$$\begin{array}{l} e_1\circ e_n = e_n\circ e_1 = e_2\circ e_n = e_n\circ e_2 = e_2\circ e_{n-1} =
e_{n-1}\circ e_3 = 0 \\[1mm]
e_n\circ e_{n-1}=e_{n-1}\circ e_n = e_n\circ e_n =0 \ \ \mbox{and}
\ \ \alpha_5=\alpha_7=\beta_3=\beta_4=0.\end{array}$$ Thus, we
obtain the following multiplications
$$\begin{array}{l} e_i\circ e_j= C_{i+j-1}^j
e_{i+j}, \  \mbox{ for }  \ 2\leq i+j\leq n-2, \end{array}$$
$$\begin{array}{llll}
e_1\circ e_{n-1}=e_n, & e_1\circ e_n=0, & e_{n-1}\circ
e_1=\alpha_6e_n, & e_{n-1}\circ
e_{n-1}=\alpha_8e_n, \\[1mm]
e_{n-1}\circ e_n=0, & e_{n-1}\circ e_2=0, & e_n\circ e_1=0, & e_n\circ e_{n-1}=0, \\[1mm]
e_n\circ e_n=0, & e_n\circ e_2=0, & e_2\circ e_{n-1}=0, & e_2\circ
e_n=0.\end{array}$$ Applying the induction and the following
chains of equalities
$$\begin{array}{l} e_n\circ e_{k+1}=e_n\circ (e_1\circ e_k)=(e_n\circ e_1)\circ e_k-e_n\circ (e_k\circ e_1)=-ke_n\circ
e_{k+1}, \\[1mm]
e_{k+1}\circ e_n=(e_1\circ e_k)\circ e_n=e_1\circ (e_k\circ
e_n)+e_1\circ (e_n\circ e_k), \\[1mm]
e_{n-1}\circ e_{k+1}=e_{n-1}\circ (e_1\circ e_k)=(e_{n-1}\circ
e_1)\circ e_k-e_{n-1}\circ (e_k\circ e_1)=\alpha_6e_n\circ
e_k-ke_{n-1}\circ e_{k+1}, \\[1mm]
e_{k+1}\circ e_{n-1}=(e_1\circ e_k)\circ e_{n-1}=e_1\circ
(e_k\circ e_{n-1})+e_1\circ (e_{n-1}\circ e_k),\end{array}$$ one
can prove $e_i\circ e_{n-1}=e_{n-1}\circ e_i=e_i\circ e_n=e_n\circ
e_i=0$ for $2\leq i\leq n-2.$\\ Thus we obtain the following
multiplication in the algebra $A:$
$$\begin{array}{l} e_i\circ e_j= C_{i+j-1}^j
e_{i+j}, \ \mbox{ for }  \ 2\leq i+j\leq n-2, \\[1mm]
e_1\circ e_{n-1}=e_n, \ \ e_{n-1}\circ e_1=\alpha e_n, \ \
e_{n-1}\circ e_{n-1}=\beta e_n.\end{array}$$

For the study the given family of algebras on isomorphity we
consider the following change of generating basic elements
$$\begin{array}{l} e'_1=ae_1+be_{n-1}, \\[1mm]
e'_{n-1}=ce_1+de_{n-1}\end{array}$$ where $ad-bc\neq 0.$

Making the change of basis
$e'_2=a^2e_2+(ab+ab\alpha+b^2\beta)e_n,$ $e'_i=a^ie_i,$ for $3\leq
i\leq n-2,$ consider the multiplication $$e'_2\circ
e'_{n-1}=2a^2ce_3.$$ On the other hand, $e'_2\circ e'_{n-1}=0,$
hence $c=0.$

Comparing in a similar way the coefficients in decomposition of
the following multiplications
$$
e'_1\circ e'_{n-1}, \ \ e'_{n-1}\circ e'_1, \ \ e'_{n-1}\circ
e'_{n-1}
$$
in bases $\{e_1, e_2, \ldots, e_n\}$ and $\{e'_1, e'_2, \ldots,
e'_n\},$ we obtain
$$d(a+b\beta)\neq 0, \,\ a\alpha+b\beta=\alpha'(a+b ), \,\
\beta'=\frac{d}{a+b\beta}\beta .$$ It should be noted validity of
the equality
$$\alpha'-1=\frac{a}{a+b\beta}(a-1).$$
If $\beta=0,$ then $\beta'=0$ and we obtain the one-parametric
family of algebras $KF_n^1.$ If $\beta\neq 0,$ then putting
$d=\frac{a+b\beta}{\beta},$ we obtain $\beta'=1.$ In the case
$\alpha=1$ we have $\alpha'=1,$ i.e. the algebra $A$ is isomorphic
to the algebra $KF_n^2.$ But if $\alpha\neq 1,$ putting
$b=-\frac{a\alpha}{\beta},$ we get $\alpha=0$ and the algebra
$KF_n^3.$ \\ \textbf{Case 1.2.} Let  $\beta_1=0$. Put
$e'_1=ae_1+be_{n-1}.$ Then
$$e'_2=e'_1\circ
e'_1=(a^2+ab\alpha_5+b^2\alpha_7)e_2+(ab(1+\alpha_6)+b^2\alpha_8)e_n.$$
Consider the multiplication $$e'_1\circ
e'_2=[(\beta_2+(1+\alpha_6)\beta_3)a^2b+(\alpha_5\beta_2+\alpha_8\beta_3+
(1+\alpha_6)\beta_4)ab^2+(\alpha_7\beta_2+\alpha_8\beta_4)b^3]e_3.$$
If at least one of the expressions  $ \beta_2+(1+ \alpha_6)
\beta_3,$ $\alpha_5 \beta_2+ \alpha_8 \beta_3+(1+ \alpha_6)
\beta_4,$ and $\alpha_7 \beta_2+ \alpha_8 \beta_4$ is different
from zero, we have case 1.1. That is why we suppose
 $$\begin{array}{c} \beta_2+(1+ \alpha_6)
\beta_3=0, \\[1mm]
\alpha_5 \beta_2+ \alpha_8 \beta_3+(1+ \alpha_6) \beta_4=0,
\\[1mm]
\alpha_7 \beta_2+ \alpha_8 \beta_4=0.\end{array}$$\\ \textbf{Case
1.2.1.} Let  $\beta_3\neq 0.$ Then by scale of basis we obtain
$\beta_3=1.$ The equalities
$$Z(e_1,e_1,e_1)=Z(e_1,e_1,e_{n-1})=Z(e_1,e_{n-1},e_1)=Z(e_2,e_1,e_1)=0$$
deduce
$$e_2\circ e_1=e_2\circ e_2=0, \ \ e_2\circ e_{n-1}=e_n\circ
e_1=(1+\alpha_6)e_3.$$ Since $\beta_2+(1+\alpha_6)=0,$ we have
$e_{n-1}\circ e_2=-(1+\alpha_6)e_3.$ Hence $$(e_{n-1}\circ
e_1)\circ e_1=-2(1+\alpha_6)e_3.$$ On the other hand,
$(e_{n-1}\circ e_1)\circ e_1=\alpha_6(1+\alpha_6)e_3,$ therefore
\begin{equation} \label{e4}
(1+\alpha_6)(2+\alpha_6)=0,
\end{equation}
One can easily obtain $$e_i\circ e_1=((i-1)+\alpha_6)e_{i+1}, \ \
 e_i\circ e_2=0, \ \ 3\leq i\leq n-3.$$
On the other hand, $$e_3\circ
e_2=\frac{1}{2}(2+\alpha_6)(3+\alpha_6)e_5, \ \ e_4\circ
e_2=\frac{1}{2}(3+\alpha_6)(4+\alpha_6)e_6.$$ Thus, we obtain
\begin{equation} \label{e5}
(2+\alpha_6)(3+\alpha_6)=0, \ \ (3+\alpha_6)(4+\alpha_6)=0
\end{equation}
It should be noted that for $n\geq 8$ the (\ref{e4}) and
(\ref{e5}) lead to a contradiction with existence of an algebra in
this case. \\ \textbf{Case 1.2.2.} Let  $\beta_3=0$. Then
$$\alpha_8= \beta_2=0, \ \ \alpha_6=-1, \ \ \beta_4\neq 0.$$
Note that taking as $e'_1$ the expression $ae_1+be_{n-1},$ we get
into case 1.2.1. In fact, at $e'_1=ae_1+be_{n-1}$ we have
$$e'_n=e'_1\circ e'_{n-1}=ae_n+b\alpha_7e_2,$$ what follows
$e'_1\circ e'_n=ab\beta_4e_3.$ If $ab\neq 0,$ then we can assume
$e'_1\circ e'_n=e'_3,$ where $e'_3=ab\beta_4e_3$ and we are in
conditions of case 1.2.1. \\ \textbf{Case 2.} Let $\alpha_4=0$ and
$\alpha_6\neq 0.$ Then we can assume $e_{n-1}\circ e_1=e_n.$ Hence
$$\begin{array}{ll} e_1\circ e_1=e_2, & e_1\circ e_{n-1}=\alpha_3e_2, \\  e_{n-1}\circ e_1=e_n, & e_{n-1}\circ
e_{n-1}=\alpha_7e_2+\alpha_8e_n.\end{array}$$ Consider the
following change of generating elements of the basis $\{e_1, e_2,
\ldots, e_n\}$ in the form:
$$\begin{array}{l}
e'_1=ae_1+be_{n-1}, \\
e'_{n-1}=ce_1+de_{n-1}
\end{array}$$ where $ad-bc\neq 0.$ Then in the new basis we have from multiplication $e'_1\circ e'_1=e'_2$
that
$$e'_2=(a^2+ab\alpha_3+b^2\alpha_7)e_2+(ab+b^2\alpha_8)e_n.$$
Moreover, we have $$e'_1\circ
e'_{n-1}=(ac+ad\alpha_3+bd\alpha_7)e_2+(bc+bd\alpha_8)e_n.$$ Note
that we can choose the numbers $a, b, c, d$ such that they satisfy
the conditions:
\begin{equation}\label{e6}a^2+ab\alpha_3+b^2\alpha_7\neq 0, \,\
bc+bd\alpha_8\neq 0.
\end{equation}
If
$\frac{a^2+ab\alpha_3+b^2\alpha_7}{ac+ad\alpha_3+bd\alpha_7}\neq
\frac{ab+b^2\alpha_8}{bc+bd\alpha_8},$ then denoting
$$e'_n=(ac+ad\alpha_3+bd\alpha_7)e_2+(bc+bd\alpha_8)e_n$$ we get into case 1.
Let $\frac{a^2+ab\alpha_3+b^2\alpha_7}{ac+ad\alpha_3+bd\alpha_7}=
\frac{ab+b^2\alpha_8}{bc+bd\alpha_8}$ for any values of $a, b, c,
d$ satisfying (\ref{e6}) at fixed $\alpha_3,$ $\alpha_7,$ and
$\alpha_8.$ Then
$$(a^2+ab\alpha_3+b^2\alpha_7)(bc+bd\alpha_8)-(ac+ad\alpha_3+bd\alpha_7)(ab+b^2\alpha_8)=0.$$
$\Rightarrow  a(ad-bc)(\alpha_3-\alpha_8)+b(bc-ad)\alpha_7=0
\Rightarrow \alpha_8=\alpha_3$ and $\alpha_7=0.$  Taking into
account that $\alpha_8=\alpha_3$ and $\alpha_7=0,$ we obtain
$e'_2=(a^2+ab\alpha_3)e_2+(ab+b^2\alpha_3)e_n \ \ \Rightarrow$
$a\neq -b\alpha_3.$

Express in the basis $\{e'_1, e'_2, \ldots, e'_n\}$ the following
multiplications:
$$\begin{array}{l}e'_1\circ
e'_{n-1}=\displaystyle \frac{c+d\alpha_3}{a+b\alpha_3}e'_2, \\[1mm]
$$e'_{n-1}\circ e'_1=(ac+bc\alpha_3)e_2+(ad+bd\alpha_3)e_n=e'_n,
\\[1mm]
e'_{n-1}\circ e'_{n-1}=\displaystyle
\frac{c+d\alpha_3}{a+b\alpha_3}e'_n.\end{array}$$ Choosing
$c=-d\alpha_3$ $(d\neq 0)$ $\Rightarrow$ $e'_1\circ
e'_{n-1}=e'_{n-1}\circ e'_{n-1}=0,$ we obtain
$$\begin{array}{l}
e_1\circ e_1=e_2, \ \ e_{n-1}\circ e_1=e_n, \\[1mm]
e_1\circ e_2=\beta_1e_3, \ \ e_{n-1}\circ e_2=\beta_2e_3, \ \
e_1\circ e_n=\beta_3e_3, \ \ e_{n-1}\circ
e_n=\beta_4e_3.\end{array}$$ \\ \textbf{Case 2.1.} Let
$\beta_1\neq 0.$ Then by scale of basis we obtain $e_1\circ
e_i=e_{i+1}$ for $1\leq i\leq n-3,$ what allows easily to obtain
validity of the equality
$$e_i\circ e_j= C_{i+j-1}^j
e_{i+j}, \  \mbox{ for }  \ 2 \leq i+j \leq n-2.$$ The equalities
$$\begin{array}{l}
Z(e_1,e_1,e_1)=Z(e_1,e_{n-1},e_1)=Z(e_1,e_1,e_{n-1})=Z(e_{n-1},e_1,e_1)= \\
=Z(e_1,e_n,e_1)=Z(e_{n-1},e_{n-1},e_1)=Z(e_{n-1},e_1,e_n)=0,\end{array}$$
yield $$\beta_2=\beta_3=\beta_4=0, \,\ e_2\circ e_1=2e_3, \,\
e_2\circ e_{n-1}=e_n\circ e_n=e_n\circ e_1=0.$$ Sum up
multiplications being in the presence:
$$\begin{array}{l} e_i\circ e_j= C_{i+j-1}^j
e_{i+j}, \ \ \mbox{for} \ \ 2 \leq i+j \leq n-2, \ \
e_{n-1}\circ e_1=e_n, \\
e_1\circ e_{n-1}=e_1\circ e_n=e_2\circ e_{n-1}=e_2\circ
e_n=e_{n-1}\circ e_2=0, \\
e_{n-1}\circ e_{n-1}=e_{n-1}\circ e_n=e_n\circ e_1=e_n\circ e_2
=e_n\circ e_{n-1}=e_n\circ e_n=0.\end{array}$$ In a similar way as
in case 1.1 one can prove $$ e_i\circ e_{n-1}=e_{n-1}\circ
e_i=e_i\circ e_n=e_n\circ e_i=0 \ \mbox{ at } \ \ 2\leq i\leq
n-2.$$ In the conclusion of this case we have that $A$ is
isomorphic to $KF_n^4.$ \\ \textbf{Case 2.2.} Let $\beta_1=0.$
Then $e_1\circ e_2=0.$ Set $e'_1=ae_1+be_{n-1}.$ Then the
condition $e'_1\circ e'_1=e'_2,$ yield $e'_2=a^2e_2+abe_n$.
Consider the multiplication $$e'_1\circ
e'_2=ab(a(\beta_2+\beta_3)+b\beta_4)e_3.$$ If either
$\beta_2+\beta_3$ or $\beta_4$ doesn't equal zero, we get into
case 2.1. Consider the case of $\beta_2+\beta_3=0$ and
$\beta_4=0.$ Then we have $\beta_2\neq 0$ (since otherwise, i.e.
at $\beta_i=0,$ $1\leq i\leq 4$ we have $dimA_3=0$). Putting
$e'_3=\beta_2e_3,$ one can assume $\beta_2=1.$
$Z(e_1,e_{n-1},e_1)=0$ $\Rightarrow$ $e_3=0,$ what contradicts to
the existence condition of an algebra in this case. \\
\textbf{Case 3.} Let $\alpha_4=\alpha_6=0.$ Then  $\alpha_8\neq 0$
and putting $e'_1=ae_1+be_{n-1},$ where $a$ and $b$ are such that
$$ab\neq 0 \ \ \mbox{and} \ \
det \left(
\begin{array}{cc}
a^2+ab\alpha_3+ab\alpha_5 & b^2 \\
\alpha_3 & b\alpha_8 \\
\end{array}
\right)\neq 0,$$ we reduce this case to the case 1.

It should be noted that we consider a fortiori non isomorphic
cases.
\end{proof}
\begin{thm} Any five-dimensional naturally graded quasi-filiform
Zinbiel algebra of type $A_{(2)}$ is isomorphic to one of the
following pairwise non isomorphic algebras:
\begin{description}
    \item[$KF_5^1:$] $\left\{\begin{array}{l} e_1\circ e_1=e_2, \ \ e_4\circ e_4=e_5,
    \ \  e_1\circ e_2=e_3, \ \ e_4\circ
e_5=e_3, \\[1mm]  e_2\circ e_1=2e_3, \ \ e_5\circ
e_4=2e_3;\end{array}\right.$
    \item[$KF_5^2:$] $\left\{\begin{array}{l}
e_1\circ e_1=e_2, \ \ e_1\circ e_2=e_3, \ \ e_1\circ e_4=e_5, \ \
e_1\circ e_5=e_3,\\[1mm] e_2\circ e_1=2e_3, \ \ e_4\circ e_1=-e_2,
\ \ e_4\circ e_2=-e_3, \ \ e_4\circ e_5=-e_3;
\end{array}\right.$
    \item[$KF_5^3(\beta):$]$
\left\{\begin{array}{l} e_1\circ e_1=e_2, \ \ e_1\circ e_2=e_3, \
\ e_1\circ e_4=e_5, \ \ e_1\circ e_5=\beta e_3,
\\[1mm]
e_2\circ e_1=2e_3, \ \ e_2\circ e_4=(\beta-1)e_3, \ \ e_4\circ
e_1=-e_2, \ \ e_4\circ e_4=-e_5,\\[1mm] e_4\circ e_2=-e_3,\ \ e_4\circ
e_5=-\beta e_3, \ \  e_5\circ e_1=(\beta-1)e_3, \ \ e_5\circ
e_4=-2\beta e_3, \ \ \beta \in C;
\end{array}\right.$
    \item[$KF_5^4:$] $\left\{\begin{array}{l}e_1\circ e_1=e_2, \ \ e_1\circ e_2=e_3, \ \
    e_1\circ e_4=e_5, \ \ e_2\circ
e_1=2e_3,\\[1mm] e_4\circ e_4=-e_5;\end{array}\right.$
    \item[$KF_5^5:$] $\left\{\begin{array}{l}e_1\circ e_4=e_2, \ \ e_1\circ e_2=e_3, \ \
    e_1\circ e_5=-e_3, \ \ e_4\circ
e_1=e_5,\\[1mm] e_4\circ e_2=e_3, \ \ e_4\circ e_5=-e_3;\end{array}\right.$

    \item[$KF_5^6:$] $e_1\circ e_4=e_2, \ \ e_1\circ e_2=e_3, \ \
    e_1\circ e_5=-e_3, \ \ e_4\circ
e_1=e_5;$
   \item[$KF_5^7:$] $\left\{\begin{array}{l} e_1\circ e_1=e_2, \ \ e_1\circ
e_2=e_3, \ \ e_1\circ e_4=e_5, \ \ e_1\circ e_5=e_3, \\[1mm]
e_2\circ e_1=2e_3, \ \ e_2\circ e_4=e_3, \ \ e_5\circ e_1=e_3;
\end{array}\right.$
    \item[$KF_5^8:$] $e_1\circ e_1=e_2, \ \ e_1\circ e_2=e_3, \ \ e_1\circ e_4=e_5, \ \ e_2\circ
e_1=2e_3;$
    \item[$KF_5^9:$] $\left\{\begin{array}{l}e_1\circ e_1=e_2, \ \ e_1\circ e_2=e_3,
    \ \ e_1\circ e_4=-e_5, \ \ e_2\circ
e_1=2e_3,\\[1mm] e_4\circ e_1=e_5, \ \ e_4\circ e_5=e_3;\end{array}\right.$
    \item[$KF_5^{10}:$] $e_1\circ e_1=e_2, \ \ e_1\circ e_4=-e_5, \ \ e_4\circ e_1=e_5, \ \ e_4\circ
e_5=e_3;$
    \item[$KF_5^{11}:$] $\left\{\begin{array}{l}e_1\circ e_1=e_2, \ \ e_1\circ e_2=e_3, \ \ e_1\circ e_4=-e_5,
    \ \ e_1\circ e_5=e_3,\\[1mm] e_2\circ e_1=2e_3, \ \ e_4\circ e_1=e_5;\end{array}\right.$
    \item[$KF_5^{12}:$] $\left\{\begin{array}{l}e_1\circ e_1=e_2, \ \ e_1\circ e_2=e_3,
    \ \ e_1\circ e_4=-e_5, \ \ e_2\circ
e_1=2e_3,\\[1mm] e_4\circ e_1=e_5;\end{array}\right.$
    \item[$KF_5^{13}:$] $e_1\circ e_1=e_2, \ \ e_1\circ e_4=-e_5, \ \ e_1\circ e_5=e_3, \ \ e_4\circ
e_1=e_5;$
    \item[$KF_5^{14}(\alpha):$] $\left\{\begin{array}{l}e_1\circ e_1=e_2,\ \
     e_1\circ e_2=e_3, \ \ e_1\circ e_4=\alpha e_5,\ \
e_2\circ e_1=2e_3,\\[1mm] e_4\circ e_1=e_5, \ \ \alpha \in C\backslash
\{-1\};\end{array}\right.$
    \item[$KF_5^{15}(\alpha):$] $\left\{\begin{array}{l} e_1\circ e_1=e_2,\ \
     e_1\circ e_4=\alpha e_5, \ \ e_1\circ
e_5=\frac {2\alpha}{1+\alpha}e_3,\ \  e_2\circ e_4 =2\alpha
e_3,\\[1mm]
e_4\circ e_1=e_5, \ \ e_4\circ e_2=e_3,\ \ e_5\circ e_1=2e_3,\ \
\alpha \in C\backslash \{-1, -\frac{1}{2} \}; \end{array}\right.$
    \item[$KF_5^{16}:$] $\left\{\begin{array}{l} e_1\circ e_1=e_2, \ \ e_1\circ
    e_2=e_3,\ \
e_1\circ e_4=-\frac {1}{2}e_5, \ \ e_1\circ e_5=-2e_3,\\[1mm] e_2\circ
e_1=2e_3, \ \ e_2\circ e_4=-e_3,\ \ e_4\circ e_1=e_5,\ \  e_4\circ
e_2=e_3,\\[1mm] e_5\circ e_1= 2e_3;\end{array}\right.$
    \item[$KF_5^{17}:$] $\left\{\begin{array}{l} e_1\circ e_1=e_2, \ \ e_1\circ e_4=-\frac
{1}{2}e_5, \ \ e_1\circ e_5=-2e_3, \ \ e_2\circ e_4 =-e_3, \\[1mm]
e_4\circ e_1=e_5, \ \ e_4\circ e_2=e_3, \ \ e_5\circ e_1= 2e_3.
\end{array}\right.$
\end{description}
\end{thm}
\begin{proof} Let $A$ be an algebra satisfying the conditions of
the theorem and let $\{e_1, e_2, e_3, e_4, e_5\}$ be a basis of
the algebra satisfying the natural gradating,
$$A_1=\langle e_1, e_4\rangle, \ \  A_2=\langle e_2, e_5\rangle, \
\ A_3=\langle e_3\rangle.$$ Write the multiplication of the basic
elements in the form
$$\begin{array}{l}
e_1\circ e_1=\alpha_1e_2 +\alpha_2e_5, \ \ e_1\circ
e_4=\alpha_3e_2+\alpha_4e_5, \\
e_4\circ e_1=\alpha_5e_2+\alpha_6e_5, \ \ e_4\circ
e_4=\alpha_7e_2+
\alpha_8e_5, \\
e_1\circ e_2=\beta_1e_3, \ \ e_1\circ e_5=\beta_2e_3, \ \ e_4\circ
e_2=\beta_3e_3, \ \ e_4\circ e_5=\beta_4e_3,\end{array}$$ where
$(\beta_1, \beta_2, \beta_3, \beta_4)\neq (0, 0, 0, 0).$ \\It is
easy to see that the linear span $\langle e_3\rangle$ is an ideal
of $A.$ Consider now the quotient algebra $A/I=\{\overline{e}_1,
\overline{e}_4, \overline{e}_2, \overline{e}_5 \}.$ It is a
four-dimensional Zinbiel algebra for which conditions $dim
(A/I)^2=2$ and $dim(A/I)^3=0$ hold. Using classification of
four-dimensional Zinbiel algebras according to Theorem
\ref{th2.3}, we conclude that $A/I$ is isomorphic to the following
pairwise non isomorphic algebras:
\begin{description}
    \item[$M_1:$] $e_1\circ e_1=e_2, \ \ e_4\circ e_4=e_5;$

    \item[$M_2:$] $e_1\circ e_1=e_2, \ \ e_1\circ e_4=e_5, \ \ e_4\circ e_1=-e_5,\ \ e_4\circ
e_4=e_2-2e_5;$
    \item[$M_3:$] $e_1\circ e_1=e_2, \ \ e_1\circ e_4=e_5, \ \ e_4\circ
    e_1=-e_2;$
    \item[$M_4(\alpha):$] $e_1\circ e_1=e_2, \ \ e_1\circ e_4=e_5, \ \ e_4\circ e_1=
\alpha e_2, \ \ e_4\circ e_4=-e_5, \ \ \alpha \in C;$
    \item[$M_5:$] $e_1\circ e_4=e_2, \ \ e_4\circ e_1=e_5;$
    \item[$M_6:$] $e_1\circ e_1=e_2, \ \ e_1\circ e_4=e_5;$
    \item[$M_7(\alpha):$] $e_1\circ e_1=e_2, \ \ e_1\circ e_4=\alpha e_5, \ \ e_4\circ
e_1=e_5 \ \ \alpha\in C.$
\end{description}
Hence we can get the values for structural constants $\alpha_i$ of
the algebra $A$, namely by equating the values of $\alpha_i$ to
the corresponding ones in the algebra from the above list.
Applying standard classification methods in each of the seven
cases complete the proof of the theorem.
\end{proof}
\newpage
\begin{thm} Any six-dimensional naturally graded quasi-filiform
Zinbiel algebra of type $A_{(2)}$ is isomorphic to one of the
following pairwise non isomorphic algebras:
\begin{description}
    \item[$KF_6^1:$] $\left\{\begin{array}{l} e_i\circ e_j=C_{i+j-1}^{j}e_{i+j},
     \ \ 2\leq i+ j\leq 4,\\[1mm]
     e_1\circ e_5=e_6, \ \ e_5\circ e_1=e_6;\end{array}\right.$
    \item[$KF_6^2:$] $\left\{\begin{array}{l} e_i\circ e_j=C_{i+j-1}^{j}e_{i+j},\ \
     2\leq i+j\leq 4,\\[1mm]
     e_1\circ e_5=e_6, \ \ e_5\circ e_1=e_6, \ \ e_5\circ e_5=e_6;\end{array}\right.$
    \item[$KF_6^3:$] $\left\{\begin{array}{l}e_i\circ e_j=C_{i+j-1}^{j}e_{i+j},
    \ \ 2\leq i+ j\leq 4,\\[1mm]
    e_1\circ e_5=e_6,\ \  e_5\circ e_5=e_6;\end{array}\right.$
    \item[$KF_6^4:$] $\left\{\begin{array}{l}e_i\circ e_j=C_{i+j-1}^{j}e_{i+j}, \ \
    2\leq i+j\leq 4,\\[1mm]
    e_5\circ e_1=e_6;\end{array}\right.$
    \item[$KF_6^5:$] $\left\{\begin{array}{l}e_1\circ e_1=e_2, \ \
    e_1\circ e_3=e_4, \ \ e_1\circ e_5=e_6, \ \ e_1\circ e_6=e_3,\\[1mm] e_2\circ
    e_5=-e_3, \ \ e_5\circ e_1=-3e_2-2e_6, \ \ e_5\circ e_2=e_3,\ \  e_5\circ
    e_5=2e_2+e_6,\\[1mm]
    e_5\circ e_6=-2e_3, \ \
    e_6\circ e_1=-e_3, \ \ e_6\circ e_5=e_3;\end{array}\right.$
    \item[$KF_6^6:$] $\left\{\begin{array}{l}e_1\circ e_1=e_2, \ \ e_1\circ e_3=e_4,
    \ \ e_1\circ e_5=e_6,\ \ e_1\circ e_6=e_3,\\[1mm] e_2\circ
    e_5=-e_3, \ \ e_5\circ e_1=-2e_6, \ \ e_5\circ e_2=e_3, \ \ e_6\circ e_1=-e_3;\end{array}\right.$
    \item[$KF_6^7:$] $\left\{\begin{array}{l}e_1\circ e_1=e_2, \ \ e_1\circ e_3=e_4,
    \ \ e_1\circ e_5=e_6,\ \ e_1\circ e_6=e_3,\\[1mm]
    e_3\circ e_1=e_4, \ \ e_5\circ e_1=-e_6.\end{array}\right.$
\end{description}
\end{thm}
\begin{proof} Similar to the cases 1.1 and 2.1 of Theorem
\ref{th4} we can get the existence (pairwise non isomorphic) of
the following algebras
$$KF_6^1, \ \ KF_6^2, \ \ KF_6^3, \ \ KF_6^4.$$
Analogously to the case 1.2 we have for the algebra $A$ with the
basis $\{e_1, e_2, e_3, e_4, e_5, e_6\}$
$$\begin{array}{c}
e_1\circ e_1=e_2, \ \ e_1\circ e_5=e_6, \ \ e_1\circ e_2=0, \ \
e_1\circ
e_6=e_3, \ \ e_1\circ e_3=e_4, \\
e_5\circ e_1=\alpha_5e_2+\alpha_6e_6, \ \ e_5\circ
e_5=\alpha_7e_2+\alpha_8e_6, \ \ e_5\circ e_2=-(1+\alpha_6)e_3, \\
e_5\circ e_6=\beta_4e_3, \ \ e_5\circ e_3=\gamma e_4\end{array}$$
where for parameters $\alpha_5, \alpha_6, \alpha_7, \alpha_8,
\beta_4, \gamma,$ the relations $$\begin{array}{c} -
5(1+\alpha_6)+\alpha_8+(1+\alpha_6)\beta_4=0, \\[1mm]
-7(1+\alpha_6)+\alpha_8\beta_4=0\end{array}$$ hold. Consideration
of the identity (\ref{e1}) for the basic elements reduces to the
following restrictions:
\begin{equation} \label{e7}
\left\{\begin{array}{ll} \alpha_8(2+\alpha_6)=0, \\[1mm]
\alpha_8(\beta_4+2\alpha_8)=0, \\[1mm]
(1+\alpha_6)(2+\alpha_6)=0, \\[1mm]
(1+\alpha_6)(\beta_4+2\alpha_8)=0, \\[1mm]
-\alpha_7(1+\alpha_6)+\alpha_8\beta_4=0, \\[1mm]
\beta_4(2+\alpha_6)=(2+\alpha_6)\gamma, \\[1mm]
\alpha_5(1+\alpha_6)+2\alpha_6\alpha_8=-\alpha_8, \\[1mm]
\beta_4(\beta_4+2\alpha_8)=(\beta_4+2\alpha_8)\gamma, \\[1mm]
-\alpha_5(1+\alpha_6)+\alpha_8+(1+\alpha_6)\beta_4=0, \\[1mm]
\alpha_5(2+\alpha_6)+\alpha_6(\beta_4+2\alpha_8)=(2+\alpha_6)\gamma,
\\[1mm]
\alpha_7(2+\alpha_6)+\alpha_8(\beta_4+2\alpha_8)=(\beta_4+2\alpha_8)\gamma.
\end{array} \right.
\end{equation} \\ \textbf{Case 1.} Let $\alpha_6=-2.$ Then (\ref{e7}) becomes to the
form
$$\left\{\begin{array}{ll}\beta_4=-2\alpha_8, \\[1mm]
\alpha_5=-3\alpha_8, \\[1mm]
\alpha_7=2\alpha_8^2.\end{array}\right.$$ Multiplication in the
algebra in this case has the form
$$\begin{array}{llll} e_1\circ e_1=e_2, & e_1\circ e_5= e_6, & e_1\circ e_6=e_3,  & e_1\circ
e_3=e_4, \\[1mm]
e_5\circ e_1=-3\alpha e_2-2e_6, & e_5\circ e_5=2\alpha^2e_2+\alpha
e_6, & e_5\circ e_2=e_3, & e_5\circ e_6=-2\alpha e_3, \\[1mm]
e_5\circ e_3=\gamma e_4, & e_2\circ e_5= -e_3, & e_6\circ
e_1=-e_3, & e_6\circ e_5=\alpha e_3,\end{array}$$ If
$\gamma+\alpha\neq 0,$ then substituting
$$\begin{array}{lll}
e'_1=(\gamma+\alpha)e_1, & e'_2=(\gamma+\alpha)^2e_2, &
e'_3=(\gamma+\alpha)^2e_3, \\[1mm]
e'_4=(\gamma+\alpha)^3e_4, & e'_5=-\gamma e_1+e_5, &
e'_6=(\gamma+\alpha)(-\gamma e_2+e_6),
\end{array}$$ we can assume $\gamma=0$ and $\alpha=1,$ i.e. we
obtain the algebra  $KF_6^5.$ But if  $\gamma+\alpha=0,$ then
substituting
$$e'_1=e_1, \ \ e'_2=e_2, \ \ e'_3=e_3, \ \ e'_4=e_4, \ \ e'_5=-\gamma e_1+e_5, \ \ e'_6=-\gamma
e_2+e_6,$$ we obtain $\gamma=0$ and $\alpha=0$, hence we have the
algebra $KF_6^6.$ \\ \textbf{Case 2.} Let $\alpha_6=-1.$ Then (7)
become the form
$$\left\{\begin{array}{ll} \alpha_8=0, \\[1mm]
\alpha_5=2\gamma, \\[1mm]
\alpha_7=\gamma^2, \\[1mm]
\beta_4 =\gamma.\end{array}\right.$$ Assume obtained table of
multiplication:
$$\begin{array}{llll} e_1\circ e_1=e_2, & e_1\circ e_5=e_6, & e_1\circ e_2=0, & e_1\circ
e_6=e_3, \\[1mm]
e_1\circ e_3=e_4, & e_5\circ e_1=2\beta e_2-e_6, & e_5\circ
e_5=\beta^2e_2, & e_5\circ e_6=\beta e_3, \\[1mm]
e_5\circ e_3=\beta e_4, & e_6\circ e_6=\beta e_4, & e_3\circ
e_1=e_4, & e_3\circ e_5=\beta e_4.\end{array}$$ Making the change
of the basis: $$e'_1=e_1, \ \ e'_2=e_2, \ \ e'_3=e_3, \ \
e'_4=e_4, \ \ e'_5=-\beta e_1+e_5, \ \ e'_6=-\beta e_2+e_6,$$ we
obtain the algebra $KF_6^7.$

Since consideration of cases 1.1, 1.2, and 2.1 was chosen such
that algebras satisfying these various cases were non isomorphic,
the sets of algebras $\{KF_6^1, KF_6^2, KF_6^3\},$ $\{KF_6^4 \},$
$\{KF_6^5, KF_6^6, KF_6^7\}$ are pairwise not disjoint (up to
isomorphism). Pairwise non isomorphity of algebras $KF_6^1,$
$KF_6^2,$ $KF_6^3$ follows from Theorem \ref{th4}. By virtue of
$dimR(KF_6^5)=dimR(KF_6^6)=1,$ and $dimR(KF_6^7)=2,$ the algebra
$KF_6^7$ is not isomorphic to algebras $KF_6^i$ $(i=5,6).$ Non
isomorphity of algebras $KF_6^5$ and $KF_6^6$ can be easily
checked by consideration of general change of the basis. Thus, we
obtain pairwise non isomorphic algebras $KF_6^i,$ $1\leq i\leq 7.$
\end{proof}
The following theorem can be proved in the same manner.

\begin{thm} Any seven-dimensional naturally graded quasi-filiform
Zinbiel algebra of type $A_{(2)}$ is isomorphic to one of the
following pairwise non isomorphic algebras:

\begin{description}
    \item[$KF_7^1:$] $\left\{\begin{array}{l}
    e_i\circ e_j=C_{i+j-1}^{j}e_{i+j}, \ \ 2\leq i+j\leq 5,\\[1mm]
    e_1\circ e_6=e_7, \ \ e_6\circ e_1=\alpha e_6;\end{array}\right.$
    \item[$KF_7^2:$] $\left\{\begin{array}{l}
    e_i\circ e_j=C_{i+j-1}^{j}e_{i+j}, \ \ 2\leq i+ j\leq 5,\\[1mm]
    e_1\circ e_6=e_7,\ \  e_6\circ e_1=e_7, \ \ e_6\circ e_6=e_7;\end{array}\right.
   $
    \item[$KF_7^3:$] $\left\{\begin{array}{l}
    e_i\circ e_j=C_{i+j-1}^{j}e_{i+j}, \ \ 2\leq i+j\leq 5,\\[1mm]
    e_1\circ e_6=e_7, \ \ e_6\circ e_6=e_7;\end{array}\right.$
    \item[$KF_7^4:$] $\left\{\begin{array}{l}e_i\circ e_j=C_{i+j-1}^{j}e_{i+j},
    \ \ 2\leq i+j\leq 5,\\[1mm]
    e_6\circ e_1=e_7;\end{array}\right.$
    \item[$KF_7^5:$] $\left\{\begin{array}{l} e_1\circ e_1=e_2, \ \
    e_1\circ e_3=e_4, \ \ e_1\circ e_4=e_5,\ \
    e_1\circ e_6=e_7,\\[1mm]
    e_1\circ e_7=e_3,\ \  e_2\circ e_6=-e_3, \ \ e_3\circ e_1=-e_4,\ \ e_6\circ
    e_1=-2e_7,\\[1mm]
    e_6\circ e_2=e_3, \ \ e_7\circ e_1= -e_3;\end{array}\right.$
    \item[$KF_7^6:$] $\left\{\begin{array}{l} e_1\circ e_1=e_2,\ \  e_1\circ e_3=e_4,
    \ \ e_1\circ e_4=e_5,\ \ e_1\circ e_6=e_7,\\[1mm]
    e_1\circ e_7=e_3,\ \  e_2\circ e_6=-e_3, \ \ e_3\circ e_1=-e_4, \ \ e_6\circ
    e_1=-2e_7,\\[1mm]
    e_6\circ e_2=e_3,\ \ e_6\circ e_4=e_5, \ \ e_7\circ e_1=-e_3.\end{array}\right.$
\end{description}
\end{thm}

\subsection{The case $r > 2$}

The proving the remaining cases, we need the following lemmas.

\begin{lem} \label{l1} Let $A$ be a naturally graded quasi-filiform
Zinbiel algebra of type $A_{(r)}$ $(r>2).$ Then $x\circ x=0$ for
any $x\in A_1.$
\end{lem}
\begin{proof} Assume to the contrary, that is there exists $x \in
A_1$ such that $x\circ x\neq 0.$ Then we choose a basis $\{e_1,
e_2, \ldots, e_n\}$ of $A$ such that $e_1=x,$ $e_2=x\circ x,$ and
$A_1=\langle e_1, e_{n-1}\rangle,$ $A_2=\langle e_2\rangle,$
$\ldots,$ $A_r=\langle e_r, e_n\rangle,$ $A_{r+1}=\langle
e_{r+1}\rangle,$ $\ldots,$ $A_{n-2}=\langle e_{n-2}\rangle .$
\\Thus, we can assume $$\begin{array}{l} e_1\circ e_i=e_{i+1} \
\mbox{ for } \ \ 2\leq i\leq r-1, \\[1mm]
e_{n-1}\circ e_{r-1}=e_n.\end{array}$$ On the other hand, similar
to the case of a filiform Zinbiel algebra in \cite{Ada}, we obtain
$e_{n-1}\circ e_i=0$ for $2\leq i\leq r-1,$ which contradicts the
existence of an element $x$ such that $x\circ x\neq 0.$
\end{proof}

\begin{lem} \label{l2} Let $A$  be a naturally graded quasi-filiform
Zinbiel algebra of type $A_{(r)}.$ Then $r\leq 3.$
\end{lem}

\begin{proof} Assume to the contrary, i.e. $r>3.$ By Lemma \ref{l1},
$x\circ x=0$ for any $x\in A_1.$
\\Choose a basis $\{e_1, e_2, \ldots, e_n\}$ of the algebra $A$ such
that
$$\begin{array}{l} e_1\circ e_i=e_{i+1} \ \ \mbox{for} \ \ 2\leq i\leq
r-1 \\[1mm]
e_1\circ e_1= e_{n-1}\circ e_{n-1}=0, \ \ e_1\circ
e_{n-1}=-e_{n-1}\circ e_1=e_2.\end{array}$$ We get a contradiction
from the equalities
$$\begin{array}{l} e_2\circ e_1=(e_1\circ e_{n-1})\circ e_1=e_1\circ (e_{n-1}\circ e_1)
+e_1\circ (e_1\circ e_{n-1})=e_1\circ (-e_2+e_2)=0, \\[1mm]
0=(e_1\circ e_1)\circ e_2=e_1\circ (e_1\circ e_2)+e_1\circ
(e_2\circ e_1) =e_1\circ e_3=e_4,\end{array}$$ thus completing the
proof of the lemma. \end{proof}

\begin{lem} \label{l3} Let $A$  be a naturally graded quasi-filiform
Zinbiel algebra of type $A_{(3)}.$ Then $dimA\leq 7.$
\end{lem}

\begin{proof} Suppose that $dimA>7$ and let $\{e_1, e_2, \ldots, e_n\}$ be a
basis satisfying the conditions $A_1=\langle e_1, e_{n-1}\rangle,$
$A_2=\langle e_2\rangle,$ $A_3=\langle e_3, e_n\rangle,$ $A_4=
\langle e_4\rangle,$ $\dots,$ $A_{n-2}=\langle e_{n-2}\rangle,$
$$\begin{array}{l} e_1\circ e_1=e_{n-1}\circ e_{n-1}=0, \\[1mm]
e_1\circ e_{n-1}=-e_{n-1}\circ e_1=e_2, \\[1mm]
e_1\circ e_2=e_3, \ \ e_{n-1}\circ e_2=e_n.\end{array}$$ The
equalities
$$\begin{array}{l} Z(e_1,e_{n-1},e_1)=Z(e_1,e_{n-1},e_{n-1})=
Z(e_1,e_1,e_2)=Z(e_{n-1},e_{n-1},e_2)= \\[1mm]
Z(e_1,e_{n-1},e_2)=Z(e_{n-1},e_1,e_2)=Z(e_1,e_2,e_{n-1})=Z(e_{n-1},e_2,e_1)=
\\[1mm]
Z(e_1,e_2,e_1)=Z(e_{n-1},e_2,e_{n-1})=0,\end{array}$$ lead to
$$\begin{array}{l} e_2\circ e_1=e_2\circ e_{n-1}=e_1\circ e_3=e_{n-1}\circ e_n=0\\[1mm]
e_3\circ e_1=e_1\circ e_3=e_n\circ e_{n-1}=e_{n-1}\circ e_n=0,
\\[1mm]
e_2\circ e_2=e_1\circ e_n=e_3\circ e_{n-1}=-e_{n-1}\circ
e_3=-e_n\circ e_1=\gamma e_4.\end{array}$$ Note that $\gamma\neq
0.$ Otherwise $A_4=0,$ i.e. $dimA\leq 5.$ \\Without loss of
generality we can assume that $\gamma=1$ and
$$\begin{array}{l} e_1\circ e_i=e_{i+1} \ \ \mbox{for} \ \
4\leq i\leq n-2, \\[1mm]
e_1\circ e_n=e_4.\end{array}$$ Using the identity (\ref{e1}) for
elements $\{e_1, e_n, e_1\},$ we obtain $e_4\circ e_1=0.$ \\
On the other hand, the equalities $$(e_1\circ e_1)\circ
e_4=e_1\circ (e_1\circ e_4)+e_1\circ (e_4\circ e_1)$$ imply
$0=e_1\circ e_5=e_6,$ i.e. we arrive at a contradiction, which
completes the proof of the lemma.
\end{proof} From Lemma \ref{l3} one can easily derive the following
corollaries.

\begin{cor} Any five-dimensional naturally graded quasi-filiform
Zinbiel algebra of type $A_{(3)}$ is isomorphic to the algebra
$$e_1\circ e_2=e_3, \ \ e_2\circ e_1=-e_3, \ \ e_1\circ
e_3=e_4, \ \ e_2\circ e_3=e_5.$$ \end{cor}
\begin{cor} Any six-dimensional naturally graded quasi-filiform
Zinbiel algebra of type $A_{(3)}$ is isomorphic to the algebra
$$\begin{array}{l} e_1\circ e_2=e_3, \ \ e_1\circ
e_3=e_4, \ \ e_1\circ e_5=e_6, \ \ e_2\circ e_1=-e_3, \ \ e_2\circ
e_3=e_5,
\\[1mm]
e_2\circ e_4=-e_6, \ \ e_3\circ e_3=e_6, \ \ e_4\circ e_2=e_6,\ \
e_5\circ e_1=-e_6.\end{array}$$
\end{cor}
\begin{cor} Any seven-dimensional naturally graded quasi-filiform
Zinbiel algebra of type $A_{(3)}$ is isomorphic to the algebra
$$\begin{array}{llll} e_1\circ e_2=e_3,& e_1\circ e_3=e_4, & e_1\circ e_5=e_6,
&  e_1\circ e_6=e_7,\\[1mm]
e_2\circ e_1=-e_3, & e_2\circ e_3=e_5, & e_2\circ e_4=-e_6,&
e_3\circ e_3=e_6,\\[1mm] e_4\circ e_2=e_6, & e_4\circ e_3=2e_7, & e_5\circ
e_1=-e_6. & \end{array}$$ \end{cor} Thus, we obtain the
classification of complex naturally graded quasi-filiform Zinbiel
algebras of an arbitrary dimension. In fact, the results of this
paper complete the classification of $n$-dimensional nilpotent
naturally graded algebras $A$ satisfying the condition
$A^{n-2}\neq 0.$


\begin{thebibliography}{8}

\bibitem{Ada} Adashev J.Q., Omirov B.A. and Khudoyberdiyev A.Kh.
On some nilpotent classes of Zinbiel algebras and their
applications. Third International Conference on Research and
Education in Mathematics. 2007. Malaysia. pp. 45-47.

\bibitem{Dzhu1} Dzhumadil'daev A.S., Identities for multiplications derived by
Leibniz and Zinbiel multiplications. Abstracts of short
communications of International conference "Operator algebras and
quantum theory of probability" (2005), Tashkent, pp. 76-77.

\bibitem{Dzhu2} Dzhumadil'daev A.S. and Tulenbaev K.M., Nilpotency of Zinbiel
algebras. J. Dyn. Control. Syst., vol. 11(2), 2005, pp. 195-213.

\bibitem{Gunz1} Gunzburg V. and Kapranov M.  Koszul duality for operads,
Duke Math. J. vol. 76, 1994, pp. 203-273.

\bibitem{Lod1} Loday J.-L., Cup product for Leibniz cohomology and dual
Leibniz algebras. Math Scand., vol. 77, 1995, pp. 189-196.

\bibitem{Lod2} Loday J.-L. and Pirashvili T., Universal enveloping algebras of
Leibniz algebras and (co)homology. Math.Ann. vol. 296, 1993, pp.
139-158.

\bibitem{O2} Omirov B.A., On derivations of filiform Leibniz algebras.
Math. Notes, v. 77(5), 2005,  pp. 733-742.

\bibitem{O3} Omirov B.A., Classification of two-dimensional complex Zinbiel algebras,
Uzbek. Mat. Zh., vol. 2, 2002, pp. 55-59.

\bibitem{Ver}
Vergne M. Cohomologie des alg\`{e}bres de Lie nilpotentes.
Application $\grave{a}$ l'\'{e}tude de la vari\'{e}t\'{e} des
alg\`{e}bres de Lie nilpotentes. Bull. Soc. Math. France, v. 98,
1970, pp. 81 - 116.
\end{thebibliography}
\end{document}